\documentclass{gSTA2e}

\usepackage{subfigure}
\pdfoutput=1
\theoremstyle{plain}
\newtheorem{theorem}{Theorem}[section]

\theoremstyle{remark}

\theoremstyle{definition}

\begin{document}


\articletype{RESEARCH ARTICLE}

\newcommand{\sbin}[2]{{\textstyle{ {#1}\choose {#2}}}}
\newcommand{\jbin}[2]{{{#1}\choose {#2}}}


\title{ \v{S}id\'ak-type tests for the two-sample problem based on precedence and exceedance statistics}

\author{Eugenia Stoimenova$^{\rm a}$$^{\ast}$\thanks{$^\ast$Corresponding author. Email: jeni@math.bas.bg
\vspace{6pt}} and N. Balakrishnan$^{\rm b}$\\\vspace{6pt}  $^{\rm a}${\em{Institute of Information and Communication Technologies and Institute of Mathematics and Informatics, Bulgarian Academy of Sciences, Acad. G.Bontchev str., block 25A, 1113 Sofia, Bulgaria}}; $^{\rm b}${\em{Department of Mathematics and Statistics, McMaster University, Hamilton, ON L8S 4K1, Canada}}\\\vspace{6pt}\received{\today } }

\maketitle

\begin{abstract}
This paper deals with a class of nonparametric two-sample tests for ordered alternatives. The test statistics proposed are based on the number of observations from one sample that precede or exceed a threshold specified by the other sample, and they  are extensions of \v{S}id\'ak's test. We derive their exact null distributions and also discuss a large-sample approximation. We then study their power properties exactly against the Lehmann alternative
and make some comparative comments. Finally, we present an example to illustrate the proposed tests.\bigskip

\begin{classcode}62G10; 62E15 \end{classcode}\bigskip

\end{abstract}

\section{Introduction}

Suppose $X$ and $Y$ are random variables with absolutely continuous univariate distributions $F$ and $G$, respectively. For testing the hypothesis $H_0: \  F(x)= G(x)$ against the alternative
\begin{equation} \label{eq:HA}
  H_A: \  F(x)> G(x),
\end{equation}
there are simple tests based on available precedences and exceedances. One can count the number of observations in the $Y$-sample above all observations in the $X$-sample, or the number of observations in the $X$-sample below all those in the $Y$-sample.


As suggested by Tukey \cite{tukey:59}, one or both of these statistics might be used to test $H_0$ against $H_A$ in (\ref{eq:HA}). The test based on the sum of these two quantities is mentioned as the earliest work of \v{S}id\'ak on nonparametric statistics; see {\cite{Seidler:00}}. The null distribution of this test statistic was  studied by \v{S}id\'ak and Vondr\'a\v{c}ek \cite{sidak:57} and tables of critical values were produced by these authors. A slight modification of the test statistic based on the sum became popular as Tukey's Quick Test (see \cite{neave:66} and  \cite{gans:81}). It basically leads to the same critical regions as \v{S}id\'ak's test. H\'ajek and \v{S}id\'ak \cite{hajek:67} found that the same statistic also leads to locally most powerful rank tests for testing $H_0$ against a one-sided shift in the location parameter if the underlying distribution is uniform. They discussed some other test statistics based on exceeding observations, such as Haga's test \citep{haga:60} and $E$-test also discussed in  \cite{sidak:77}. In all these tests, the counts were with respect to the extreme order statistics from one or both samples.

\bigskip

The extreme sample values may get inflated by possible outliers, which may adversely affect the performance of these test statistics. For this reason, we may want to reduce their influence by defining thresholds above the smallest and below the largest observed values in the samples. Let $X_1,\ldots,X_{m}$ and $Y_1,\ldots,Y_{n}$ be independent random samples from continuous distributions $F$ and $G$, respectively. Thresholds based on the $(r+1)$-th order statistic from the $Y$-sample and $(m-s)$-th order statistic from the $X$-sample define the exceedance and precedence statistics of the form
\begin{equation}\label{exst}
\begin{array}{rclcl}
A_s= & \mbox{the number of}& Y\mbox{-observations} & \mbox{larger than} & X_{(m-s)}, \\
B_r= & \mbox{the number of}& X\mbox{-observations} & \mbox{smaller than} & Y_{(1+r)},
\end{array}
\end{equation}
where $0\le s <m$ and $0\le r <n$.

In this paper, we propose a family of rank statistics for the two-sample problem in which the test statistic is a sum of $A_s$ and $B_r$ for appropriate choices of $s$ and $r$. It includes  \v{S}id\'ak's test as a special case.


\bigskip

Tests based on the number of precedences ($B_r$) were recommended by  \cite{nelson:63} for life-testing since a location shift can be effectively detected before all the data have been collected. They can be successfully applied to the general two-sample problem stated above. Some basic references on precedence tests include \cite{shorack:67}, \cite{katzen:85}, \cite{katzen:89}, \cite{chakr-lann:97}, \cite{vanderLaan:01}. There are many extensions of precedence tests; see \cite{balakr:05}, \cite{bairamov:06}, \cite{bairtanil:08}, \cite{anna:08}, \cite{balatri:08}. For more details on these developments, one may refer to Ng and Balakrishnan \cite{balakr:06}. Recently, a family of tests based on the minimum of $A_s$ and $B_r$ has been studied in \cite{jspi:11}.

\bigskip

The rest of this paper is organized as follows.  In Section \ref{thetest}, we introduce the new test statistics. In Section \ref{thenull}, we derive the exact null distributions of these test statistics and suggest some approximations for large samples. In Section \ref{thealt}, we derive the exact distributions of the test statistics under the Lehmann alternative and study the power functions of the tests against this alternative.
In Section  \ref{discu}, we compare the powers of the proposed tests with other known tests based on exceedances, and also present an illustrative example.  Proofs of the theorems are relegated to the Appendix.

\section{The proposed test statistics}\label{thetest}

To test  $H_0$ versus $H_A$ in (\ref{eq:HA}), we propose the test statistic
\begin{equation}\label{vrstat}
    V_{\rho} = A_s + B_r,
\end{equation}
where the threshold statistics $X_{(m-s)}$ and $Y_{(1+r)}$ are determined as $s=[\rho m]$ and $r=[\rho n]$ for some $0\le \rho <1$, with $[ \cdot ]$ denoting  the integer part. Various values of $\rho$ yield a family of test statistics which we refer to as \v{S}id\'ak-type tests. Reasonable values of $\rho$ are between 0 and 1/2. For $\rho=0$, it is  equivalent to \v{S}id\'ak's statistic \cite{sidak:57}. $\rho>0$ determines a part from the ordered samples that are skipped before the threshold is specified. Its role will be discussed in more detail later in Section \ref{infl}.

Evidently, large values of $V_{\rho}$ lead to the rejection of $H_0$ in favor of the stochastically ordered alternative in $H_{A}$. It is reasonable to select $\rho$ to be small since we want to reduce the possible influence of a small number of potential outliers.

For equal sample sizes, the parameters $s$ and $r$, specifying the threshold positions, are equal and in this case the contiguous order statistics determine the family of test statistics. For simplicity, let us denote the family of test statistics in this case by $V_{r}= A_r + B_r$ with $r=0, 1,2, \ldots$.

The following example is useful for an illustration of the proposed $V_{\rho}$-test statistic. The data is a subset of a data on breakdown times (in minutes) of an insulating fluid that is subjected to high voltage stress presented in \cite{nelson:82}. Take $X$- and $Y$-samples to be Samples 3 and 6 from \citep[][p. 462]{nelson:82}, respectively.

\smallskip

\noindent {\bf Example 1.} Ten units each of group $X$ and group $Y$ were placed simultaneously on a life-testing experiment, and their lifetimes (in minutes) were observed and are as presented in Table \ref{nel:p462}.
\begin{table}
 \tbl{Lifetimes of two samples of an insulating fluid.}
{\begin{tabular}{@{}l*{10}{c}}
\toprule
Group & \multicolumn{10}{c}{Lifetimes}\\
\colrule
$X$ & 0.49 & 0.64 & 0.82 & 0.93 & 1.08 & 1.99 & 2.06 & 2.15 & 2.57 & 4.75\\
$Y$ & 1.34 & 1.49 & 1.56 & 2.10 & 2.12 & 3.83 & 3.97 & 5.13 & 7.21 & 8.71\\
   \botrule
\end{tabular}}
\label{nel:p462}
\end{table}

In this case, we have $m=n=10$. Let $r=0, 1$ and 2 and take consecutively the threshold values to be the pairs $(Y_{(1)}, X_{(10)})$, $(Y_{(2)},X_{(9)})$, and $(Y_{(3)},X_{(8)})$. We find the corresponding precedence and exceedance statistics as presented in Table \ref{nel:p462a}.
\begin{table}
 \tbl{Computation of $V_{\rho}$-statistic.}
{\begin{tabular}{@{}lll}
\toprule
$r$ & $Y$-threshold & Precedences \\
\colrule
0 & $Y_{(1)} = 1.34$, & $B_0=5$, \\
1 & $Y_{(2)} = 1.49$, & $B_1=5$, \\
2 & $Y_{(3)} = 1.56$, & $B_2=5$, \\
   \botrule
\end{tabular}
\hspace{1cm}
 \begin{tabular}{@{}lll}
\toprule
$r$ & $X$-threshold & Exceedances \\
\colrule
0  & $X_{(10)} = 4.75$, & $A_0=3$, \\
1  & $X_{(9)}  = 2.57$, & $A_1=5$, \\
2  & $X_{(8)}  = 2.15$, & $A_2=5$. \\
   \botrule
\end{tabular}}
\label{nel:p462a}
\end{table}
With these, the first three \v{S}id\'ak-type test statistics are readily found to be $V_0 = 8$, $V_1 = 10$ and $V_2 = 10$.

The tests from the family (\ref{vrstat}) have advantage to some of the other rank tests in the case when a small number of outliers are expected to be present in the data.


\section{Null distribution}\label{thenull}

In this section, we derive the exact null distribution of the \v{S}id\'ak-type test statistic defined in (\ref{vrstat}), provide some tables of critical values for some selected small sample sizes, and finally suggest some approximation for large sample sizes.

\subsection{Exact distribution}

Given the joint distribution of $A_{s}$ and $B_{r}$ under the null hypothesis $H_0$,  the cumulative distribution function of $V_{\rho}$-statistic, for $0 \le z \le m+n$, is given by
\begin{equation}\label{distr0}
P(V_{\rho}\le z|F=G) =  \sum_{i=0}^{z}  \sum_{k=0}^{z-i} P(A_{s}=k,B_{r}=i|F=G).
\end{equation}

\begin{theorem}\label{th:joint0}

For any $0\le s < m$ and $0\le r < n$, the joint probability mass function of $A_{s}$ and $B_{r}$, under $H_0: F(x) = G(x)$, is given by
\begin{eqnarray*}
  P(A_{s}=k,B_{r}=i) &=& \frac{\jbin{s+k}{s}\jbin{r+i}{r}}{\jbin{m+n}{n}} \jbin{m+n-s-r-i-k-2}{n-r-k-1}, \\
    & & \mbox{for } 0\le i \le m-s-1, \ \mbox{ and } 0\le k \le n-r-1, \\
    &=&\frac{\jbin{m+n-r-i-1}{n-r-1}\jbin{m+n-s-k-1}{m-s-1}}{\jbin{m+n}{n}}
 \jbin{k+i-m-n+s+r}{k-n+r}, \\
 & & \mbox{for }  m-s \le i \le m, \ \mbox{and } n-r \le k \le n,\\
&=& 0,   \   \mbox{otherwise. }
\end{eqnarray*}

\end{theorem}

The proof of this theorem is presented in the Appendix.

\bigskip

To compute the cumulative distribution function of $V_{\rho}$-statistic under $H_0$, we just substitute for $P(A_{s}=k,B_{r}=i|F=G)$ from Theorem \ref{th:joint0}  into (\ref{distr0}).

The cumulative distribution function in (\ref{distr0}) is thus distribution-free. However, it does not take a simpler expression. For the simplest case when $\rho=0$, \v{S}id\'ak and Vondr\'a\v{c}ek \cite{sidak:57} presented the formula
\[
P(V_0 \le z) = \left\{ \jbin{m+n-z}{n} + \sum_{j=0}^{z-1} \jbin{m+n-z-1}{m-j}\right\}/\jbin{m+n}{n}.
\]

\subsection{Critical values }

Using the exact null distribution in (\ref{distr0}), we can determine the critical region of the test statistic $V_{\rho}$ for a pre-fixed level of significance $\alpha$. Under the alternative hypothesis that $Y$ is stochastically larger than $X$ as in (\ref{eq:HA}), we expect the $X$-observations to take on most of the smaller ranks. Hence, $H_0$ is rejected if $V_{\rho} \ge c$, where critical value $c$ is determined as the minimal $c$ such that $P(V_{\rho} \ge c| H_0) \le \alpha$.

For small sample sizes, the expression in (\ref{distr0}) is easy to compute
\footnote{These and further calculations have been carried out on a PC computer by using the statistical package R. The code can be provided by the corresponding author upon request.}. Table \ref{crit25} presents the critical values $c$ of the $V_{r}$-tests for the choices of the sample sizes $m=n=6,\ldots, 25$ for $\alpha = 0.05$, where the index $r$ corresponds to the threshold statistics $Y_{(1+r)}$ and $X_{(m-r)}$.

\begin{table}
\tbl{Critical values of \v{S}id\'ak-type tests for $m=n= 6(1)25$ and different choices of $r$ at 5\% level of significance}
{\begin{tabular}{c*{20}{c@{\hspace{5pt}}}}
\toprule
$r/n$ &    6 &  7 &  8&   9&  10&  11&  12 & 13 & 14 & 15&  16 & 17 & 18&  19 & 20 & 21&  22 & 23 & 24&  25  \\
\colrule
0   &    5 &  5 &  5&   5&   5&   5&   5 &  5 &  5 &  5&   5 &  6 &  6&   6 &  6 &  6&   6 &  6 &  6&   6  \\
1   &    8 &  8 &  8&   8&   8&   9&   9 &  9 &  9 &  9&   9 &  9 &  9&   9 &  9 &  9&   9 &  9 &  9&   9  \\
2   &   10 & 11 & 12&  11&  11&  11&  11 & 12 & 12 & 12&  12 & 12 & 12&  12 & 12 & 12&  12 & 12 & 12&  12  \\
3   &    * &  * & 13&  14&  15&  14&  14 & 14 & 14 & 14&  14 & 14 & 14&  14 & 15 & 14&  15 & 15 & 15&  15  \\
4   &    * &  * &  *&   *&  16&  16&  17 & 18 & 17 & 17&  17 & 17 & 17&  17 & 17 & 17&  17 & 17 & 17&  17  \\
5   &    * &  * &  *&   *&   *&  18&  18 & 19 & 20 & 20&  20 & 20 & 20&  20 & 20 & 20&  20 & 20 & 20&  20   \\
6   &    * &  * &  *&   *&   *&   *&   * & 20 & 21 & 22&  22 & 23 & 22&  22 & 22 & 22&  22 & 23 & 23&  22   \\
7   &    * &  * &  *&   *&   *&   *&   * &  * &  * & 23&  24 & 24 & 25&  25 & 26 & 25&  25 & 25 & 25&  25   \\
8   &    * &  * &  *&   *&   *&   *&   * &  * &  * &  *&  25 & 26 & 26&  27 & 27 & 28&  28 & 27 & 28&  28   \\
9   &    * &  * &  *&   *&   *&   *&   * &  * &  * &  *&   * &  * &  *&  28 & 29 & 29&  30 & 30 & 31&  30   \\
10  &    * &  * &  *&   *&   *&   *&   * &  * &  * &  *&   * &  * &  *&   * &  * & 31&  31 & 32 & 32&  33   \\
\botrule
\end{tabular}}\label{crit25}

\end{table} 

Note that due to the discreteness of the distributions of non-randomized test statistics based on ranks, the significance levels of the different $V_{r}$-tests are not the same. In order to achieve the same level of significance for all tests under study, we use the randomized test procedure described below. This allows us to make meaningful and more reasonable comparison of their power performance.

In order to achieve the prescribed $\alpha$ for all tests, for each realization (labeled $i$-th, say) of two random samples, we calculate the probability $P_i$ of rejecting $H_0$ as follows:
\begin{equation}\label{rand-prp}
    P_i = \left\{ \begin{array}{cl}
                    1, & \quad if \ V_{\rho}\ge c \\
\dfrac{\alpha - \alpha_1}{\alpha_2 - \alpha_1}, &  \quad if\ V_{\rho}=c-1 \\
                    0, &   \quad otherwise,
                  \end{array} \right.
\end{equation}
where $c$ is a possible critical value of the statistic $V_{\rho}$ such that $P(V_{\rho}\ge c)=\alpha_1$, $P(V_{\rho}\ge c-1)=\alpha_2$, with $\alpha_1 < \alpha < \alpha_2$.
For example, Table \ref{crit2} presents the exact levels of significance $\alpha_1$ and $\alpha_2$ (when $\alpha= 5\%$) of the $V_{\rho}$-tests for  $m=40$ and $n=20(4)40$.

\begin{table}
 \tbl{Critical values for $m=40$ and $n=20(4)40$ and different choices of $s$ and $r$ at 5\% level of significance}
{ \begin{tabular}{cccccccc|cccccccc}\toprule
$\rho$ &$m$ & $s$ & $n$ & $r$ &c.v. & $\alpha_1$&  $\alpha_2$
&$\rho$ &$m$ & $s$ & $n$ & $r$ &c.v. & $\alpha_1$&  $\alpha_2$\\
\colrule
  0 &  40  &  0  & 20  &  0   & 7 & 0.043 &  0.068  &   0.15 &40 &  6 & 20 &  3 & 20  & 0.044  &    0.060  \\
  0 &  40  &  0  & 24  &  0   & 7 & 0.030 &  0.050  &   0.15 &40 &  6 & 24 &  3 & 19  & 0.042  &    0.058   \\
  0 &  40  &  0  & 28  &  0   & 6 & 0.034 &  0.069  &   0.15 &40 &  6 & 28 &  4 & 20  & 0.048  &    0.067   \\
  0 &  40  &  0  & 32  &  0   & 6 & 0.034 &  0.061  &   0.15 &40 &  6 & 32 &  4 & 20  & 0.044  &    0.061   \\
  0 &  40  &  0  & 36  &  0   & 6 & 0.032 &  0.058  &   0.15 &40 &  6 & 36 &  5 & 22  & 0.037  &    0.051   \\
  0 &  40  &  0  & 40  &  0   & 6 & 0.032 &  0.058  &   0.15 &40 &  6 & 40 &  6 & 23  & 0.041  &    0.056   \\
\colrule
 0.05 &40  & 2 & 20  & 1 & 12 &  0.041&  0.059  &      0.2 & 40 & 8 &  20 &  4 & 24 &  0.044  &    0.057   \\
 0.05 &40  & 2 & 24  & 1 & 11 &  0.043&  0.064  &      0.2 & 40 & 8 &  24 &  4 & 23 &  0.041  &   0.055    \\
 0.05 &40  & 2 & 28  & 1 & 11 &  0.034&  0.053  &      0.2 & 40 & 8 &  28 &  5 & 24 &  0.043  &   0.058    \\
 0.05 &40  & 2 & 32  & 1 & 10 &  0.048&  0.075  &      0.2 & 40 & 8 &  32 &  6 & 25 &  0.047  &   0.064    \\
 0.05 &40  & 2 & 36  & 1 & 10 &  0.048&  0.075  &      0.2 & 40 & 8 &  36 &  7 & 27 &  0.039  &   0.052    \\
 0.05 &40  & 2 & 40  & 2 & 12 &  0.041&  0.062  &      0.2 & 40 & 8 &  40 &  8 & 28 &  0.043  &   0.056    \\
\colrule
 0.1 & 40 &  4 & 20  & 2 & 16 & 0.044  &   0.063  &    0.25 &40 & 10 & 20 &  5 & 28  & 0.042  &  0.053   \\
 0.1 & 40 &  4 & 24  & 2 & 15 & 0.044  &   0.062  &    0.25 &40 & 10 & 24 &  6 & 28  & 0.049  &  0.065   \\
 0.1 & 40 &  4 & 28  & 2 & 15 & 0.037  &   0.053  &    0.25 &40 & 10 & 28 &  7 & 29  & 0.048  &  0.062   \\
 0.1 & 40 &  4 & 32  & 3 & 16 & 0.046  &   0.066  &    0.25 &40 & 10 & 32 &  8 & 30  & 0.049  &  0.064   \\
 0.1 & 40 &  4 & 36  & 3 & 16 & 0.044  &   0.064  &    0.25 &40 & 10 & 36 &  9 & 32  & 0.040  &  0.051   \\
 0.1 & 40 &  4 & 40  & 4 & 18 & 0.036  &   0.052  &    0.25 &40 & 10 & 40 & 10 & 33  & 0.043  &  0.055   \\
   \botrule
\end{tabular}}\label{crit2}
\end{table}

For the $V_0$-test, \v{S}id\'ak and Vondr\'a\v{c}ek \cite{sidak:57} presented tables up to $m=26$ and $n=26$ at the 5\% and 1\% levels of significance.

For small values of $m$ and $n$, the critical values and the exact significance probabilities of the $V_{\rho}$-test can be computed without any difficulty, as done in Table \ref{crit25}. However, for large sample sizes, this would require a heavy computational effort and time. For this reason, we present below some large-sample approximations for the null distributions of  $V_{\rho}$-statistics.

\subsection{Large-sample approximation}

For the $V_0$-test, \v{S}id\'ak and Vondr\'a\v{c}ek \cite{sidak:57} presented tables of approximate critical values for 5\% and 1\% levels of significance. As $m/n \to 1$, the right tail probability of the test statistic is asymptotically equivalent to $\dfrac{c+2}{2^{c+1}}$, where $c$ is the corresponding critical value. For $\rho>0$, the tail approximation by negative binomial distribution turns out to be reasonable.
Since the  majority of the probability mass of $V_{\rho}$ is in the lower tail of the distribution, we calculate the upper tail as  $1-\sum_{i+k\le z} P(A_{s}=k,B_{r}=i)$.

\begin{theorem}\label{large-appr}
As $m,n\to \infty$ and $m/n\to 1$,
\[
P(V_{\rho}\le z) =  \sum_{i+k\le z}P(V_{\rho}= z)
\sim \sum_{i=0}^{z} \jbin{2(s+1)+i-1}{i} 2^{-2(s+1)}2^{-i},
\]
where $s=[\rho m]$.
\end{theorem}

The approximating probability is then the value given by the c.d.f. of a negative binomial random variable with parameters $2(s+1)$ and 1/2 (Figure \ref{nbappr}, left). The proof of this theorem is presented in the Appendix.


\begin{figure}
\begin{center}
\resizebox*{14cm}{!}{\includegraphics{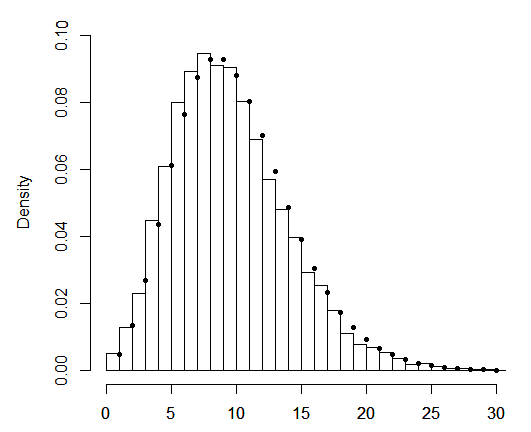} \includegraphics{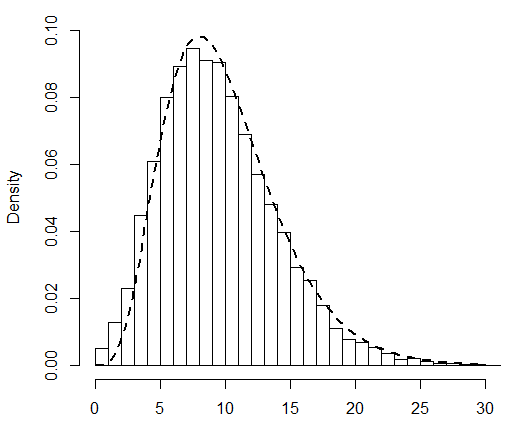}}%
\caption{Negative binomial and chi-square approximations of $V_{\rho}$ distribution for $m=400$, $\rho=0.01$.}%
\label{nbappr}
\end{center}
\end{figure}

The chi-square approximation (Figure \ref{nbappr}, right) is also quite reasonable in the practical range of sample sizes (between 25 to 100) as long as $n$ does not differ too much from $m$. In Table \ref{gam-appr2}, we provide an example of the exact significance probabilities for the $V_{\rho}$-statistics (close to 5\% level) for the choices of the sample size $m=n=40$ and 100. It is given by a chi-square distribution with degrees of freedom $[\rho m]+1$.

\begin{table}
\tbl{Values of  $P(\chi^2_{\nu} > c)$   (near 5\% critical values)}
{\begin{tabular}{cccc|cccc}
\toprule
$m$  & $\rho$ &   c.v. &  $\chi^2$-approx.  &  $m$  & $\rho$ &   c.v.  & $\chi^2$-approx. \\
\colrule
40 & 0    & 6   & 0.0497 &100&  0       &  6  & 0.0489\\
   & 0.05 & 12  & 0.0571 &   &  0.05    & 20  & 0.0649 \\
   & 0.1  & 18  & 0.0496 &   &  0.1     & 33  & 0.0587 \\
   & 0.15 & 23  & 0.0535 &   &  0.15    & 45  & 0.0585 \\
   & 0.2  & 28  & 0.0538 &   &  0.2     & 58  & 0.0475 \\
   & 0.25 & 33  & 0.0529 &   &  0.25    & 69  & 0.0513 \\
\botrule
\end{tabular}}\label{gam-appr2}
\end{table}


\section{Distributions under alternatives}\label{thealt}

\subsection{Distribution  under Lehmann alternative}\label{slehm}

In this section, we derive an expression for the distribution of $V_{\rho}$ under the Lehmann alternative given by
\begin{equation}\label{leta}
H_{LE}: \, G(x)= 1-(1-F(x))^{1/\eta},
\end{equation}
for some $\eta>1$. When $\eta=1$, the resulting distributions satisfy the null hypothesis $H_0$, while $\eta>1$ yields various distributions in the alternative hypothesis $H_{LE}$, with larger values of $\eta$ indicating stronger attraction towards $H_A: \  F(x)\ge G(x)$;  see \cite{lehmann:53} for further discussion on this class of alternatives.

As in the derivation of the null distribution, the joint probability mass function of $A_{s}$ and $B_{r}$ under $H_{LE}$ can be used for obtaining the distribution of  $V_{\rho}$ under $H_{LE}$.

Under the hypothesis $H_{LE}$ in (\ref{leta}),  the exact cumulative distribution function of the $V_{\rho}$-statistic, for $0 \le z \le m+n$,  is given by (\ref{distr0}) with the joint distribution of $A_{s}$ and $B_{r}$ now being under $H_{LE}$, as established in the following theorem.

\begin{theorem}\label{lem-lem}

For any $0\le s < m$ and $0\le r < n$, the joint probability mass function of $A_{s}$ and $B_{r}$, under $H_{LE}$ in (\ref{leta}), is given by

\begin{eqnarray*}
  P(A_{s}=k,B_{r}=i) &=& \frac{m!n!( 1/\eta)}{r!s!(n-k-r-1)!k!} S_p S_z, \\
   &&  \qquad \qquad \mbox{for } 0\le i \le m-s-1, \ \mbox{and } 0\le k \le n-r-1, \\
    &=& \frac{m!n!\eta}{(n-r-1)!(m-s-1)!(i-m+s)!(m-i)!} S_p' S_z', \\
    && \qquad\qquad \mbox{for }  m-s \le i \le m, \ \mbox{and } n-r \le k \le n,\\
    &=& 0,  \  \mbox{otherwise},
\end{eqnarray*}
where $S_p$, $S_z$, $S_p'$ and $S_z'$ are as follows:
\begin{eqnarray*}
  S_p&=& \sum_{p=0}^r (-1)^p {r\choose p}
  \frac{\Gamma(m-i+(n-r+p)/\eta)}{\Gamma(m+(n-r+p)/\eta+1)},\label{sp1}\\
  S_z&=& \sum_{z=0}^{n-k-r-1} (-1)^z {n-k-r-1\choose z}
  \frac{\Gamma(s+(z+k)/\eta+1)}{\Gamma(m-i+(z+k)/\eta+1)},\label{ss1}\\
  S_p'&=& \sum_{p=0}^{i-m+s} (-1)^p {i-m+s\choose p}
  \frac{\Gamma(n-r+(m-i+p)\eta)}{\Gamma(k+(m-i+p)\eta+1)},\label{sp2}\\
  S_z'&=& \sum_{z=0}^{m-s-1} (-1)^z {m-s-1\choose z}
  \frac{\Gamma(k+(z+s+1)\eta)}{\Gamma(n+(z+s+1)\eta+1)}.\label{ss2}
\end{eqnarray*}

\end{theorem}

The proof of this theorem is presented in the Appendix.

Consequently, the distribution of $V_{\rho}$-statistic under $H_{LE}$ is distribution-free as well.

\subsection{Power against Lehmann alternative}\label{s-lem-pow}

Now, we demonstrate the use of the exact cumulative distribution function of $V_{\rho}$ under Lehmann alternative  as well as the Monte Carlo simulation method for the computation of the power of the $V_{\rho}$-test against this alternative. For this purpose, we generated 100,000 sets of data from $F$ and $1-(1-F(x))^{1/\eta}$, respectively, and computed the test statistic $V_{\rho}$ for each set. The power values were estimated by the rejection rates of the null hypothesis for different values of $\eta$.

To make meaningful comparison of the power values of different tests, we calculated power functions at prescribed exact level  of significance $\alpha$ as follows. First, for any $V_{\rho}$-test, we determine two values $\alpha_1$ and $\alpha_2$ such that
\[
 P(V_{\rho}\ge c)=\alpha_1 \quad
  \text{and}
  \quad P(V_{\rho}\ge c-1)=\alpha_2,
\]
where $c$ is given by $P(V_{\rho}\ge c|H_0) \le \alpha$, so that the interval $(\alpha_1,\alpha_2)$ contains the critical level, say $\alpha=0.05$. Next, we calculate the power values corresponding to the two critical values $c$ and $c-1$ as
\[
  \beta_1 = P(V_{\rho}\ge c|H_{LE})
  \quad
  \text{and}
  \quad
  \beta_2 = P(V_{\rho}\ge c-1|H_{LE}).
\]
Then, the power of the test at exact level $\alpha$ is estimated by
\[
\beta= \pi\beta_2 + (1-\pi)\beta_1,
\]
where $\pi=\dfrac{\alpha - \alpha_1}{\alpha_2 - \alpha_1}$  is the adjusting factor used in the randomization  procedure  in~\eqref{rand-prp}.

For $m=n=10$ and $\eta=2(1)7$, the power values of the  $V_{\rho}$-tests corresponding to $r=0, \ldots,4$, against the Lehmann alternative  $H_{LE}$ in (\ref{leta}), are presented in Table \ref{vr:pow10}, where the significance level is set as $\alpha=0.05$.
\begin{table}
\tbl{Power comparison of $V_{r}$-tests for $m=n=10$ at 5\% level of significance}
{\begin{tabular}{ccccccc}
\toprule
$V_{r}$-test & $\eta=2$& $\eta=3$ & $\eta=4$  & $\eta=5$& $\eta=6$ & $\eta=7$\\
\colrule
$V_0$     &  0.3212  &  0.5799 & 0.7432 & 0.8415 & 0.8969  & 0.9318   \\
$V_1$     &  0.3291  &  0.5854 & 0.7430 & 0.8370 & 0.8911 &  0.9219  \\
$V_2$     &  0.3070  &  0.5536 & 0.7133 & 0.8114 & 0.8728 &  0.9097  \\
$V_3$     &  0.2946  &  0.5384 & 0.7012 & 0.8020 & 0.8673 &  0.9064   \\
$V_4$     &  0.3211  &  0.5801 & 0.7492 & 0.8468 & 0.9021 &  0.9375   \\
\botrule
  \end{tabular}}\label{vr:pow10}
\end{table}
Similarly, for $m=n=20$, the power values of the $V_{\rho}$-tests corresponding to $r = 0,\ldots, 8$, are presented in Table \ref{vr:pow20}.

\begin{table}
\tbl{Power comparison of $V_{r}$-tests for $m=n=20$ at 5\% level of significance}
{  \begin{tabular}{ccccccc}
\toprule
$V_{r}$-test & $\eta=2$& $\eta=3$ & $\eta=4$  & $\eta=5$& $\eta=6$ & $\eta=7$\\
\colrule
$V_0$  & 0.4566   &  0.7859   & 0.9207   &  0.9705 & 0.9894 & 0.9952         \\
$V_1$  & 0.5061   &  0.8292   & 0.9436   &  0.9808 & 0.9931 & 0.9974         \\
$V_2$  & 0.5230   &  0.8379   & 0.9476   &  0.9818 & 0.9928 & 0.9969         \\
$V_3$  & 0.5182   &  0.8355   & 0.9445   &  0.9795 & 0.9918 & 0.9957         \\
$V_4$  & 0.5149   &  0.8262   & 0.9416   &  0.9774 & 0.9901 & 0.9956         \\
$V_5$  & 0.4971   &  0.8137   & 0.9323   &  0.9742 & 0.9890 & 0.9948         \\
$V_6$  & 0.4737   &  0.7934   & 0.9208   &  0.9692 & 0.9866 & 0.9936         \\
$V_7$  & 0.4499   &  0.7684   & 0.9063   &  0.9618 & 0.9826 & 0.9919         \\
$V_8$  & 0.4791   &  0.8061   & 0.9299   &  0.9737 & 0.9888 & 0.9955         \\
\botrule
  \end{tabular}}\label{vr:pow20}
\end{table}

From Tables \ref{vr:pow10} and \ref{vr:pow20}, we see that the power values of all tests increase with increasing $\eta$. The power of $V_0$ (original \v{S}id\'ak test) is much less than the power of the next two $V_{r}$-tests for sample size $m=10$, and much less than the power of the next four $V_{r}$-tests for sample size $m=20$. For each of the six fixed values 2 (1) 7 of $\eta$, the power increases up to the third $V_{r}$-test, showing that the $V_0$-test, based on the extremal thresholds, is less powerful than the tests based on the next extremal thresholds pairs $(Y_{(2)}, X_{(m-1)})$ and $(Y_{(3)}, X_{(m-2)})$.

\begin{table}
 \tbl{Power of $V_{\rho}$-test against $H_1:\ G= 1- (1-F)^{1/2}$ for $m=40$
 and  $n=20(4)40$ at 5\% level of significance}
{\begin{tabular}{ccccccc}
   \toprule
  &  \multicolumn{6}{c}{Second sample size ($n$) } \\
   \colrule
   proportion   ($\rho$)        & 20 & 24 & 28 & 32 & 36 & 40\\
   \colrule
 0     & 0.3472& 0.4147 & 0.4771 & 0.5275 & 0.5685 & 0.6016  \\
 0.05  & 0.4647& 0.5242 & 0.5995 & 0.6548 & 0.7010 & 0.7367  \\
 0.1   & 0.5161& 0.5784 & 0.6313 & 0.6910 & 0.7276 & 0.7708  \\
 0.15  & 0.5489& 0.5960 & 0.6618 & 0.6937 & 0.7425 & 0.7750  \\
 0.2   & 0.5579& 0.5990 & 0.6586 & 0.7066 & 0.7428 & 0.7703  \\
 0.25  & 0.5675& 0.6203 & 0.6669 & 0.7012 & 0.7334 & 0.7625  \\
   \botrule
 \end{tabular}}\label{lem-uneq}
\end{table}

 \begin{table}
 \tbl{Power of $V_{\rho}$-test against $H_1:\ G= 1- (1-F)^{1/2}$ for $m=100$
 and $n=50(10)100$ at 5\% level of significance}
{\begin{tabular}{ccccccc}
\toprule
  &  \multicolumn{6}{c}{Second sample size ($n$) } \\
\colrule
   proportion   ($\rho$)        & 50 & 60 & 70 & 80 & 90 & 100\\
\colrule
 0     & 0.5002 & 0.5805& 0.6658 & 0.7137 & 0.7611 & 0.7858 \\
 0.05  & 0.7492 & 0.8404& 0.8936 & 0.9261 & 0.9502 & 0.9646 \\
 0.1   & 0.8205 & 0.8806& 0.9305 & 0.9569 & 0.9716 & 0.9813 \\
 0.15  & 0.8296 & 0.9030& 0.9406 & 0.9618 & 0.9750 & 0.9852 \\
 0.2   & 0.8555 & 0.9176& 0.9481 & 0.9595 & 0.9781 & 0.9825 \\
 0.25  & 0.8615 & 0.9189& 0.9439 & 0.9652 & 0.9750 & 0.9809 \\
\botrule
 \end{tabular}}\label{lem-uneq2}
\end{table}

For unequal sample sizes, we compare the power functions for fixed $\eta=2$; for other values of $\eta>1$, we observed a similar behavior and so we do not present the corresponding results for conciseness. Table \ref{lem-uneq}  provides the power values for $m=40$ and  $n=20,24,28,32,36,40$.
The proportion coefficient $\rho$ specifies the six $V_{\rho}$-tests. The power functions were estimated through Monte Carlo simulations, with  100,000 simulated data sets for each case.

\begin{figure}
\begin{center}
\resizebox*{11cm}{!}{\includegraphics{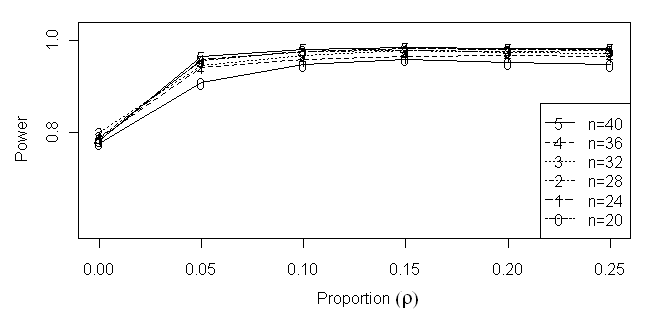}}%
\caption{Power functions of $V_{\rho}$-tests for $m=40$ and $n=20,24,28,32,36,40$ against the Lehmann alternative with $\eta=2$  at 5\% level of significance.}%
\label{lem-v-s}
\end{center}
\end{figure}

Figure \ref{lem-v-s} illustrates the gain in power of using any of the first five $V_{\rho}$-tests with $\rho>0$ instead of \v{S}id\'ak's $V_0$-test. Table \ref{lem-uneq2} provides similar results for  $m=100$ and  six  values for $n$ corresponding to six values of the proportion coefficient $\rho$.

\section{Discussion}\label{discu}

\subsection{Remark on  consistency of the test}

Since the distribution of the $V_{\rho}$-test under $H_{LE}$ is distribution-free, any particular underlying distribution $F$ can be used to prove the test consistency. Sen \cite{sen:65} has proved (in above notation) the following: For $G(x) = F(x-\theta)$ with $H_0:\, \theta=0$, the test based on $A_s$ is consistent against the set of alternatives $H_A:\, \theta>0$ for $F$ belonging to the domain of attraction for maxima of the Gumbel (type 1) family of cdf's. Consequently, the test based on the sum of $A_s$ and $B_r$ is consistent for this family of cdf's.

For our purpose, let $F$ be a Gumbel distribution, i.e.,
\[
F(x) = 1 - exp(- e^{x})
\]
for $-\infty <x< \infty$. This distribution belongs to the above mentioned family of cdf's.

The location shift alternative $G$ given by $G(x) = F(x-\theta)$ in this case is a Lehmann alternative of the form (\ref{leta}), with $1/\eta = e^{-\theta}$. Using Sen's result, we may conclude that the $V_{\rho}$-test is consistent against Lehmann alternatives.

\medskip

\subsection{Outlier-inflated distribution}\label{infl}

Suppose there are a small number (say, less than 10\%) of ``spurious'' values in the observed data set. Let us consider the following example.

\medskip

\begin{figure}
\begin{center}
\resizebox*{10cm}{!}{\includegraphics{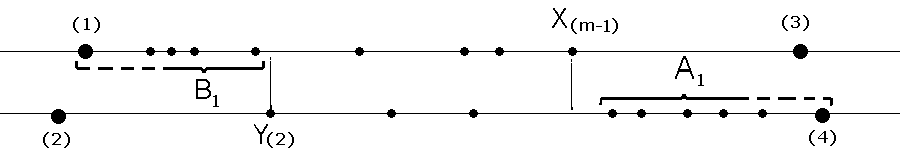}}%
\caption{Outliers present in data.}%
\label{outliers}
\end{center}
\end{figure}

\noindent {\bf Example 2.} In Figure \ref{outliers}, points labeled by (1) and (3) in the $X$-sample and (2) and (4) in the $Y$-sample lie away from the majority of observed data, i.e., they are potential outliers;  see \cite{barnett:94} for a thorough discussion on outliers.

For the hypothesis testing problem $H_0:\  F(x)=G(x)$ against $H_A: \  F(x) > G(x)$, we might want to apply some quick tests like $V_0=A_0+B_0$ or $V_1 = A_1 + B_1$. Outliers like (1) and (4) do not add much to the test statistic $V_0$; here, $B_0=0$ and $A_0=1$, while outliers like (2) and (3) inflate the thresholds and might significantly decrease $B_0$ and/or $A_0$. For this reason, the $V_{\rho}$-test with $\rho>0$ may be better since it is robust to the presence of a small number of outliers in the data.

\smallskip

\noindent {\bf Example 3.}  Let us now consider the data arising  from a contaminated distribution of the form
\[
F_{\varepsilon}  = (1-\varepsilon)F + \varepsilon F_c,
\]
where $\varepsilon$ specifies a small part of contamination with  distribution $F_c$.
Let the distribution of the second sample similarly be
\[
G_{\varepsilon}  = (1-\varepsilon)G + \varepsilon G_c.
\]

\begin{figure}
\begin{center}
\resizebox*{7cm}{!}{\includegraphics{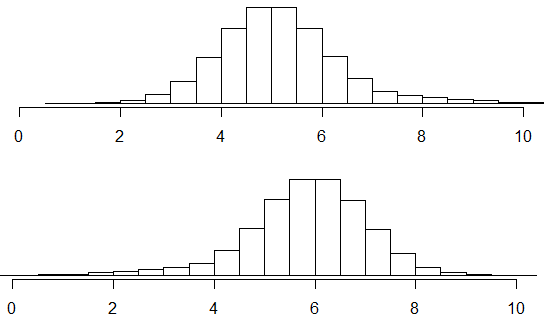}}%
\caption{Contaminated data.}%
\label{contd}
\end{center}
\end{figure}

To allow 5\% outliers in this setup, we generated samples from contaminated normal distributions as follows:
\begin{eqnarray*}
 X \sim F_{\varepsilon}  &=& 0.95\ N(5,1) + 0.05\ N(8,1), \\
  Y \sim G_{\varepsilon} &=& 0.95\ N(6,1) + 0.05\ N(3,1).
\end{eqnarray*}
The two distributions are plotted in Figure \ref{contd}.

For testing $H_0: F(x)=G(x)$ vs $H_1: F(x) > G(x)$, we use the $V_{\rho}$-tests with $\rho=0,\ 0.05,\  0.1,\ 0.15,\ 0.2,\ 0.25$. The corresponding  threshold values $s$ and $r$ are   $s=[\rho m]$ and $r=[\rho n]$, respectively. Simulating 100 observations from each distribution, we calculated the test statistics (see Table \ref{contam}).

\begin{table}
 \tbl{Test results for contaminated normal data and comparison of $V_{\rho}$-tests. (near 5\% critical values}
{\begin{tabular}{ccccccccc}
\toprule
$r=[\rho m]$ & $V_{\rho}$ & crit. value    &  \\
  \colrule
0  &    1  &  6  &                  \\
5  &    19 &  20 & Do not reject    \\
10 &    33 &  33 &                  \\
  \colrule
15 &    66 &  45 &                  \\
20 &    82 &  58 &    Reject         \\
25 &    96 &  69 &                   \\
\botrule
\end{tabular}}\label{contam}
 \end{table}

As we can expect, the first $V_{\rho}$-tests are not robust to the presence of outliers. More specifically, the $V_{\rho}$-test for $\rho=0,\ 0.05,\  0.1$ would not reject $H_0$ at 5\% level  of significance, while the $V_{\rho}$-test for $\rho=0.15,\ 0.2,\ 0.25$ suppresses the effect of outliers and do indeed reject~$H_0$.

Depending on the expected percentage of contaminated data, we could recommend to use a $V_{\rho}$-test with a suitable choice of $\rho$. Clearly, such a test will reject the null hypothesis more precisely when it is not true and its power will be similar to the power of other tests from the family. Therefore, in the case when some percentage of outliers is expected, the use of \v{S}id\'ak-type tests would be recommended.

\subsection{Comparative comments}\label{thecomp}

In this section, we discuss briefly  several nonparametric exceedance-type tests from the literature, and compare the proposed \v{S}id\'ak-type tests with these tests through an example. For more details about these tests, we refer the readers to  \cite{balakr:06}.


The classical precedence test and the maximal precedence test are useful in the case of life-testing experiments wherein data become available naturally in order of size. However, they can be used for testing $H_0: F(x)=G(x)$ against the stochastically ordered alternative as well.

\begin{itemize}
  \item For fixed $0\le r\le n$, the classical precedence test $P_r$ is simply (in terms of exceedance statistics defined by (\ref{exst})) the number of failures from the $X$-sample before the $(r+1)$-th failure from the $Y$-sample;

  \item The maximal precedence statistic $Q_r$ has been defined by  \cite{balakr:00} as the maximum number of failures occurring from the $X$-sample before the first, between the first and the second, \ldots, and between the $r$-th and $(r+1)$-th failures from the $Y$-sample;

  \item The $M_{r}$-test statistic is given by
\[
    M_{r} = \max \{ n-A_s, m-B_r \},
\]
and it was recently introduced by \cite{jspi:11}. It generalizes in some sense the $E$-test of~\cite{hajek:67};

  \item The Wilcoxon's rank-sum statistic $WR$ is known to provide a good nonparametric test for the hypothesis testing problem described above against the alternative $H_1: F(x) > G(x)$. Its test statistic is based on the sum of the ranks of observations from one of the samples obtained from the combined sample.

\end{itemize}

\begin{figure}
\begin{center}
\resizebox*{8cm}{!}{\includegraphics{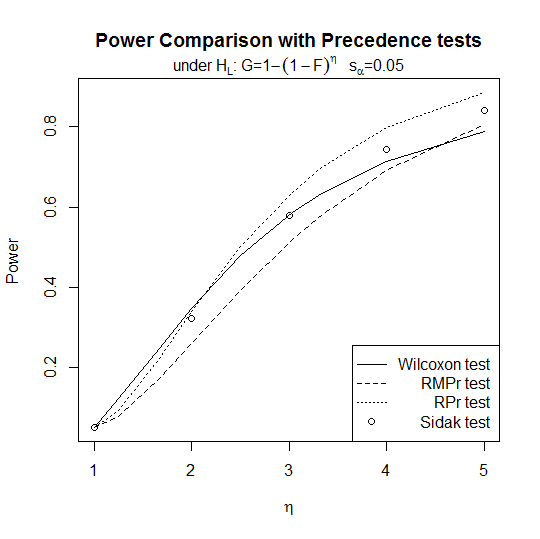}}%
\caption{Power comparison of $V_4$-test with precedence tests against Lehmann alternatives with $\eta=2$ to 5 at 5\% level of significance.}%
\label{compprec}
\end{center}
\end{figure}

The power of a $V_4$-test is compared with the power of $RPr(r)$ and $RMPr(r)$ for the case $r=4$ and  $m=n=10$. The computations here were carried for $RPr(r)$ and $RMPr(r)$ tests through $\eta = 2$ to 5 at $\alpha=0.05$. The plots are given in Figure \ref{compprec}.  Clearly, the power of the $V_r$-test is similar to the power of the two precedence-type tests and  the Wilcoxon rank-sum test. Therefore, in the case of Lehmann alternatives, the use of \v{S}id\'ak-type tests would be recommended.

In the following example, we compare the $V_{r}$-tests with the precedence test $P_r$, the maximal precedence test $Q_r$, and the $M_{r}$-test described above.

\smallskip

\noindent {\bf Example 4.} Considering the same data as in Example 1 (see Table  \ref{nel:p462}), we can carry out a nonparametric test for the hypothesis $H_0:\  F(x)=G(x)$ through the first four tests from each of the above families of tests. Table \ref{example2} provides the values of the test statistics and the corresponding p-values.

\begin{table}
 \tbl{Different test statistics and $p$-values for the insulating fluid  data}
{\begin{tabular}{ccccccccc}
\toprule
$r$ & $P_r$ & p-value &$Q_r$ & p-value & $M_r$ & p-value &$V_r$ & p-value\\
  \colrule
0 & 5 & 0.0163 & 5 &0.0163 & 7  & 0.0186  & 8  & 0.01054    \\
1 & 5 & 0.0704 & 5 &0.0325 & 5  & 0.0177  & 10 & 0.02826   \\
2 & 5 & 0.1749 & 5 &0.0487 & 5  & 0.0795  & 10 & 0.10847   \\
\botrule
\end{tabular}}\label{example2}
 \end{table}

In this example, the first five smallest $X$-values occurred before the smallest $Y$-value, and in addition, the last three largest $Y$-values occurred after the largest $X$-value. All the tests with $r=0$ perform similarly, giving evidence against $H_0$ at the usual 5\% level of significance. If $M_0$ or $V_0$ test is used, the data would provide strong evidence to reject $H_0$. However, if $P_1$-test had been used instead, it would not reject $H_0$ while the $Q_1$, $M_1$ and $V_1$ tests all would reject $H_0$. For $r \ge 2$, all tests provide similar conclusions.

\section{Acknowledgements}

The authors thank the Editor, an Associate Editor and two referees for constructive suggestions and encouragements and comments on an earlier version of the manuscript.

The work of the first author was supported by the grant I02/19 of the Bulgarian National Science Fund, while the work of the second author was supported by the Natural Sciences and Engineering Research Council of Canada through an individual discovery grant


\bibliographystyle{gSTA}
\bibliography{excee-14}


\appendices

\section{Proofs of Theorems }

\noindent{\bf Proofs of Theorems \ref{th:joint0} and \ref{lem-lem}}

Let $X_1,\ldots,X_{m}$ and $Y_1,\ldots,Y_{n}$ be two independent random samples from continuous distributions $F$ and $G$, respectively. For $0\le s <n$ and $0\le r <m$, let $A_{s}$ and $B_{r}$ be the statistics as defined in (\ref{exst}).

Consider $k$ exceedances in the $Y$-sample with respect to $X_{(m-s)}$ and $i$ precedences in the $X$-sample with respect to $Y_{(1+r)}$. First, suppose that $1\le k \le n-r-2$ and $0\le i \le n-r-1$. Event $\{A_{s}=k\}$ means that the $(m-s)$-th ordered observation from the $X$-sample is between the $(n-k)$-th and $(n-k+1)$-th ordered observations from the $Y$-sample, while event $\{B_{r}=i\}$ means that the $(1+r)$-th ordered observation from the $Y$-sample is between the $i$-th and $(i+1)$-th ordered observations from the $X$-sample. The first two cases in each of Theorems \ref{th:joint0} and \ref{lem-lem} arise according to the ordering of $Y_{(1+r)}$ and $X_{(m-s)}$:
\begin{equation}\label{th-cas}
\begin{array}{rll}
   \mbox{ If } & Y_{(1+r)}< X_{(m-s)},& \mbox{then } B_r \le m-s-1 \mbox{ and } A_s\le n-r-1, \\
\mbox{  while if } & Y_{(1+r)} > X_{(m-s)}, & \mbox{then } B_r \ge m-s \mbox{ and } A_s\ge n-r.
\end{array}
\end{equation}

Here, we derive $P(A_{s}=k,B_{r}=i)$ in the case $Y_{(1+r)} < X_{(m-s)}$ for arbitrary absolutely continuous distributions $F$ and $G$, and for the null hypothesis $F=G$ and then for the Lehmann alternative in (\ref{leta}).

Conditional on the $Y$-observations
\begin{equation}\label{cond}
Y_{(1+r)}=y_1,\  Y_{(n-k)}=y_2, \  Y_{(n-k+1)}=y_3,
\end{equation}
define the event $W_{q,t}$ on the $X$-sample as follows:
\[
W_{q,t}:= \left\{
\begin{array}{l}
    \mbox{$i$ $X$-observations preceding $y_1$}  \\
    \mbox{$t$ $X$-observations between $y_1$ and $y_2$}  \\
    \mbox{$m-i-q-t$ $X$-observations between $y_2$ and $y_3$} \\
    \mbox{$q$ $X$-observations exceeding $y_3$},
\end{array}
\right.
\]
where $0\le q\le s$ and $0\le t \le m-i-s-1$.

The probability of $W_{q,t}$ is evidently given by the multinomial probability
\begin{equation}\label{wqt}
\frac{m!}{i!t!(m-i-q-t)!q!} [F(y_1)]^i [F(y_2) - F(y_1)]^t
  [F(y_3) - F(y_2)]^{m-i-q-t} [1-F(y_3)]^q,
\end{equation}
for $y_1<y_2<y_3$. The conditional probability of  $\{A_{s}=k,B_{r}=i\}$, given (\ref{cond}), is obtained by summing (\ref{wqt}) over all $q=0, \cdots, s$ and $t=0, \cdots, m-i-s-1$. Hence, the unconditional probability of  $\{A_{s}=k,B_{r}=i\}$, with respect to the joint distribution of $Y_{(1+r)}$, $Y_{(n-k)}$ and $Y_{(n-k+1)}$, is
\begin{eqnarray}
\lefteqn{P(A_{s}=k,B_{r}=i) } \nonumber\\
&=&\sum_{q=0}^{s}\sum_{t=0}^{m-i-s-1}
 \frac{m!}{i!t!(m-i-q-t)!q!}
\int\limits_{-\infty}^{\infty} \int\limits_{y_1}^{\infty}  \int\limits_{y_2}^{\infty}
[F(y_1)]^i [F(y_2) - F(y_1)]^t\nonumber\\
& &
 \times [F(y_3) - F(y_2)]^{m-i-q-t}[1-F(y_3)]^q
g_{*} (y_1, y_2, y_3) \, dy_3 \, dy_2\, dy_1,\label{joint1}
 \end{eqnarray}
where $g_{*}$ is the joint density function of the three order statistics $Y_{(1+r)}$, $Y_{(n-k)}$ and $Y_{(n-k+1)}$, from the $Y$-sample given by (see \cite{David:03} or   \cite{arnold:08})
\begin{eqnarray}
g_{*} (y_1, y_2, y_3)&= &\frac{n!}{r!(n-k-r-2)!(k-1)!} [G(y_1)]^{r} [G(y_2) - G(y_1)]^{n-k-r-2}\nonumber\\
&& \times
 [1 - G(y_3)]^{k-1} g(y_1) g(y_2) g(y_3), \qquad \mbox{for } \ y_1<y_2<y_3,\label{jodens}
\end{eqnarray}
with $g$ being the density corresponding to $G$.

\bigskip

\noindent{\bf Proof of Theorem \ref{lem-lem}}

Under the Lehmann alternative in (\ref{leta}), the two distributions satisfy the relationships $(1-G) = (1-F)^{1/\eta}$ and $g(x)=(1/\eta)[1-F(x)]^{(1/\eta)-1}f(x)$, with $f$ and $g$ being the densities corresponding to $F$ and $G$, respectively. Substituting these in  (\ref{joint1}) and then using (\ref{jodens}), we obtain
\begin{eqnarray*}
\lefteqn{ P(A_{s}=k,B_{r}=i|H_{LE})= \sum_{q=0}^{s}\sum_{t=0}^{m-i-s-1}C
\int\limits_{-\infty}^{\infty} \int\limits_{y_1}^{\infty}  \int\limits_{y_2}^{\infty}
[F(y_1)]^i [F(y_2) - F(y_1)]^t }\\
& & \times\ [F(y_3) - F(y_2)]^{m-i-q-t}[1-F(y_3)]^q[1-(1-F(y_1))^{1/\eta}]^r [ (1-F(y_1))^{1/\eta} \nonumber\\
&& - (1-F(y_2))^{1/\eta}]^{n-k-r-2}
[1 - F(y_3)]^{(k-1)/\eta}[1 - F(y_1)]^{(1/\eta)-1}\nonumber\\
&& \times\ [1 - F(y_2)]^{(1/\eta)-1}[1 - F(y_3)]^{(1/\eta)-1}
 f(y_1) f(y_2) f(y_3) \, dy_3 \, dy_2\, dy_1,\nonumber\\
 && \label{star1}
 \end{eqnarray*}
where \ $\displaystyle C=\frac{m!n!(1/\eta)^3}{i!t!(m-i-q-t)!q!r!(n-k-r-2)!(k-1)!}$.

Changing variables in the integral by $u_i=1-F(y_i)$, $i=1,2,3$, together with $du_i=-f(y_i)dy_i$, we get
\begin{eqnarray}
\lefteqn{ P(A_{s}=k,B_{r}=i|H_{LE})}\nonumber\\
& =& \sum_{q=0}^{s}\sum_{t=0}^{m-i-s-1} C
 \int\limits_{0}^{1} \int\limits_{0}^{u_1}  \int\limits_{0}^{u_2}
(1-u_1)^i (u_1-u_2)^t (u_2-u_3)^{m-i-q-t}u_3^q\nonumber\\
&&
\hspace{1cm} \times\ (1-u_1^{1/\eta})^{r}(u_1^{1/\eta} - u_2^{1/\eta})^{n-k-r-2} u_3^{(k-1)/\eta} u_1^{1/\eta-1}
u_2^{1/\eta-1} u_3^{1/\eta-1} \, du_3 \, du_2 \, du_1\nonumber\\
&& \nonumber\\
 &=&\sum_{q=0}^{s}\sum_{t=0}^{m-i-s-1} C
 \sum_{z=0}^{n-k-r-2}  (-1)^z \sbin{n-k-r-2}{z} \sum_{p=0}^{r} (-1)^p \sbin{r}{p}
 \int\limits_{0}^{1} \int\limits_{0}^{u_1}  \int\limits_{0}^{u_2}
(1-u_1)^i\nonumber\\
&&
\hspace{1.6cm}\times\ (u_1-u_2)^t (u_2 -u_3)^{m-i-q-t} u_3^q u_1^{p/\eta} u_2^{z/\eta}
u_1^{(n-k-r-2-z)/\eta}  u_3^{(k-1)/\eta}\nonumber\\
 &&\hspace{1.6cm} \times\ u_1^{1/\eta-1}  u_2^{1/\eta-1} u_3^{1/\eta-1}   \, du_3 \, du_2 \, du_1,\label{eqalt}
 \end{eqnarray}
where in the last expression we have used binomial expansions for the power terms $(u_1^{1/\eta} - u_2^{1/\eta})^{n-k-r-2}$ and $(1-u_1^{1/\eta})^{r}$.

Then the integral in (\ref{eqalt}) is simplified by the substitution $w=u_3/u_2$ with $du_3=u_2\, dw$ and further by $w=u_2/u_1$ with $du_2=u_1\, dw$, yielding
\begin{eqnarray*}
J &=& \int\limits_{0}^{1}  \int\limits_{0}^{u_1} \int\limits_{0}^{u_2}
(1-u_1)^i (u_1-u_2)^t (u_2-u_3)^{m-i-q-t} u_1^{(n-k-r-1-z+p)/\eta-1} u_2^{(z+1)/\eta-1}\\
&&  \hspace{8.5cm} \times u_3^{q+k/\eta-1}
du_3\, du_2\, du_1\nonumber\\
&=& \int\limits_{0}^{1}  \int\limits_{0}^{u_1} \int\limits_{0}^{1}
(1-u_1)^i (u_1-u_2)^t u_2^{m-i-q-t} (1-w)^{m-i-q-t}
u_1^{(n-k-r-1-z+p)/\eta-1}\nonumber\\
&&  \hspace{4.6cm} \times u_2^{(z+1)/\eta-1} u_2^{q+k/\eta-1} w^{q+k/\eta-1}u_2\,
dw\, du_2\, du_1\\
&=& B(q+k/\eta, m-i-q-t+1)\int\limits_{0}^{1} \int\limits_{0}^{u_1} (1-u_1)^{i}(u_1-u_2)^t
u_2^{m-i-t+(k+z+1)/\eta-1} \nonumber\\
&&   \hspace{7.2cm}  \times
 u_1^{(n-k-r-1-z+p)/\eta-1}    du_2\, du_1\\
&=&B(q+k/\eta, m-i-q-t+1) \int\limits_{0}^{1}  \int\limits_{0}^{1} (1-u_1)^i u_1^t  (1-w)^{t} u_1^{m-i-t+(k+z+1)/\eta-1}\\
& &  \hspace{3.6cm} \times
w^{m-i-t+(k+z+1)/\eta-1} u_1^{(n-k-r-1-z+p)/\eta-1} u_1 \, dw \, du_1\\
&=&
B(q+k/\eta, m-i-q-t+1)B(m-i-t+(k+z+1)/\eta,t+1)\\
&&   \times B(m-i+(n-r+p)/\eta,i+1),
\end{eqnarray*}
where $B(a,b)=\int_0^1 t^{a-1}(1-t)^{b-1}\,dt$ denotes the complete beta function.

Now substituting the above expression for $J$ in (\ref{eqalt}) and expressing the beta functions through gamma functions, we obtain
\begin{eqnarray*}
\lefteqn{ P(A_{s}=k,B_{r}=i|H_{LE}) =\sum_{z=0}^{n-k-r-2}
(-1)^z \jbin{n-k-r-2}{z} \sum_{p=0}^{r}(-1)^p \jbin{r}{p}}\nonumber\\
& & \times
\sum_{q=0}^{s}\sum_{t=0}^{m-i-s-1} \frac{m!n!(1/\eta)^3}{q!r!(n-k-r-2)!(k-1)!}\times\frac{\Gamma(q+k/\eta) }{\Gamma(m-i-t+k/\eta+1)}\nonumber\\
&& \times \frac{\Gamma(m-i-t+(k+z+1)/\eta)}{\Gamma(m-i+(k+z+1)/\eta+1)} \frac{\Gamma(m-i+(n-r+p)/\eta)}{\Gamma(m+(n-r+p)/\eta+1)}\\
&=& \frac{m!n!(1/\eta)^3}{r!(n-k-r-2)!(k-1)!}
 \left(\sum_{p=0}^{r}(-1)^p \jbin{r}{p} \frac{\Gamma(m-i+(n-r+p)/\eta)}{\Gamma(m+(n-r+p)/\eta+1)}\right)\\
&&  \times \left(\sum_{q=0}^{s} \frac{\Gamma(q+k/\eta) }{q!}\right)\sum_{z=0}^{n-k-r-2} (-1)^z \sbin{n-k-r-2}{z} \frac{1}{\Gamma(m-i+(k+z+1)/\eta+1)}
\\
&&  \times\left( \sum_{t=0}^{m-i-s-1}\frac{\Gamma(m-i-t+(k+z+1)/\eta)}{\Gamma(m-i-t+k/\eta+1)}\right).\label{eqsu}
\end{eqnarray*}

The sums in the above expression were simplified as follows. The sum over $t$ is
\begin{eqnarray*}
\lefteqn{Q_1 = \sum_{t=0}^{m-i-s-1}
\frac{\Gamma(m-i+(k+z+1)/\eta-t)}{\Gamma(m-i+k/\eta+1-t)}}\\
&& =- \frac{\eta}{(z+1)} \left[\frac{\Gamma(s+(k+z+1)/\eta+1)} {\Gamma(s+k/\eta +1)}
- \frac{\Gamma(m-i+(k+z+1)/\eta+1)}{\Gamma(m-i+k/\eta+1)} \right],
\end{eqnarray*}
where we used the identity
$\displaystyle \quad
\sum_{k=0}^n \frac{\Gamma(a-k)}{\Gamma(b-k)} = \frac{1}{b-a-1} \left[
\frac{\Gamma(a-n)}{\Gamma(b-n-1)} - \frac{\Gamma(a+1)}{\Gamma(b)}\right]$ \quad for $a>n$ and $b>n+1$, $b\neq a+1$ \citep[p. 492]{prudnikov:02}.

\smallskip

Next, we used the identity
$\displaystyle \quad
\sum  _{k=0}^n (-1)^k \jbin{n}{k}\frac{1}{k+p} = \frac{1}{p}\jbin{n+p}{n}^{-1}
$\quad \citep[see][p. 498]{prudnikov:02}
so that the sum over $z$ can be simplified as
\begin{eqnarray*}
 Q_2 & = & \sum_{z=0}^{n-k-r-2} (-1)^z \jbin{n-k-r-2}{z}
\frac{1}{\Gamma(m-i+(k+z+1)/\eta+1)} Q_1\\
&=& \left[  \sum_{z=0}^{n-k-r-2} (-1)^z \jbin{n-k-r-2}{z} \frac{1}{(z+1)/\eta}\right]
\frac{1}{\Gamma(m-i+k/\eta+1)} -\frac{1}{\Gamma(s+k/\eta+1)} \\
\end{eqnarray*}
\begin{eqnarray*}
&&\times\left[  \sum_{z=0}^{n-k-r-2} (-1)^z \jbin{n-k-r-2}{z}
\frac{\Gamma(s+(z+1+k)/\eta+1)}{\Gamma(m-i+(z+1+k)/\eta+1)}
\frac{1}{(z+1)/\eta}\right]
 \\
 &=&\frac{\eta}{\Gamma(m-i+k/\eta+1)} \sum_{z=0}^{n-k-r-2} (-1)^z \jbin{n-k-r-2}{z}\frac{1}{z+1}-\frac{\eta}{\Gamma(s+k/\eta+1)}\\
&& \times \frac{1}{(n-k-r-1)}
 \sum_{z=0}^{n-k-r-2} (-1)^z \jbin{n-k-r-1}{z+1}
 \frac{\Gamma(s+(z+1+k)/\eta+1)}{\Gamma(m-i+(z+1+k)/\eta+1)}\\
&=&\frac{\eta}{\Gamma(s+k/\eta+1)(n-k-r-1)}
 \sum_{z=0}^{n-k-r-1} (-1)^z \jbin{n-k-r-1}{z}\nonumber\\
&&  \times\frac{\Gamma(s+(z+k)/\eta+1)}{\Gamma(m-i+(z+k)/\eta+1)},
\end{eqnarray*}
and similarly the sum over $q$ can be simplified as
\[
Q_3 = \sum_{q=0}^{s} \frac{\Gamma(q+k/\eta)}{q!}
 = \frac{\eta}{k} \frac{\Gamma(s+k/\eta+1)}{s!},
\]
where we have used the identity
$\displaystyle \quad
\sum_{k=0}^n \frac{\Gamma(k+a)}{\Gamma(k+1)} = \frac{1}{a} \frac{\Gamma(n+a+1)}{\Gamma(n+1)}
$ \quad
for any $a>0$.

Thus, we obtain
\begin{eqnarray}
\lefteqn{ P(A_{s}=k,B_{r}=i|H_{LE})} \nonumber\\
&=& \frac{m!n! (1/\eta)^3\Gamma(s+k/\eta+1)\eta}{r!(n-k-r-2)!(k-1)!k s!}
\frac{\eta}{\Gamma(s+k/\eta+1)} \frac{1}{(n-k-r-1)} S_p S_z \nonumber\\
&=& \frac{m!n! (1/\eta)}{r!s!(n-k-r-1)!k!} S_p S_z, \label{joint3}
\end{eqnarray}
where $S_p$ and $S_z$ are given by
\begin{eqnarray*}
  S_p&=& \sum_{p=0}^r (-1)^p {r\choose p}
  \frac{\Gamma(m-i+(n-r+p)/\eta)}{\Gamma(m+(n-r+p)/\eta+1)},\\
  S_z&=& \sum_{z=0}^{n-k-r-1} (-1)^z {n-k-r-1\choose z}
  \frac{\Gamma(s+(z+k)/\eta+1)}{\Gamma(m-i+(z+k)/\eta+1)}.
\end{eqnarray*}

We have thus derived the first case of Theorem \ref{lem-lem} for $0\le i \le m-s-1$ and $2\le k \le n-r-2$. It can be easily extended for $k=0$ and $k= n-r-1$ by using the joint density of two order statistics from distribution~$G$.

In the second case of Theorem \ref{lem-lem}, when $ m-s\le i\le n$ and $ n-r\le k\le n$, the ordering of the observations can be viewed as a symmetric image of the ordering for the first case with the following switches: $F\leftrightarrow G$;\ $(r+1) \leftrightarrow (m-s)$; \ $i \leftrightarrow (n-k)$;\ $k \leftrightarrow (m-i)$.

So, if we denote $\psi(m,n,s,r,\eta,k,i)$ to be the RHS of (\ref{joint3}), i.e,
\[
\psi(m,n,s,r,\eta,k,i)=\frac{m!n!( 1/\eta)}{r!s!(n-k-r-1)!k!} S_p S_z,
\]
then for $ m-s\le i\le n$ and $ n-r\le k\le n$, we have
\begin{eqnarray}
 P(A_{s}=k,B_{r}=i|H_{LE})&=&\psi(n,m,n-r-1,m-s-1, 1/\eta,m-i,n-k) \nonumber\\
 && \hspace{-2cm}
= \displaystyle \frac{m!n!\eta}{(n-r-1)!(m-s-1)!(i-m+s)!(m-i)!} S_p' S_z',\label{joint3a}
\end{eqnarray}
where $S_p'$ and $S_z'$ are as stated in the theorem.

The last case  of the theorem follows trivially due to (\ref{th-cas}).

\bigskip

\noindent{\bf Proof of Theorem \ref{th:joint0}}

The proof of the theorem follows readily by substituting  $\eta =1$ in Theorem \ref{lem-lem}. Hence,  for $0\le i \le m-s-1$ and $0\le k \le n-r-1$, the sums  $S_p$ and $S_z$ have closed-forms and after simplification, they become
\begin{eqnarray*}
S_{p0} & = & \frac{1}{(i+1)!} \sum_{p=0}^r (-1)^p {r\choose p} {m+n-r+p \choose i+1}^{-1} \\
& = & \frac{1}{(i+1)!} \frac{i+1}{r+i+1}  {m+n \choose m+n-r-i-1}^{-1}
 =\frac{(m+n-r-i-1)!(r+i)!}{i!(m+n)!},  \\
  S_{z0} &=& \frac{1}{(m-i-s)!}\sum_{z=0}^{n-k-r-1} (-1)^z {n-k-r-1\choose z} {m-i+k+z \choose m-i-s}^{-1}  \\
& = &  
\frac{(k+s)!(m+n-s-r-i-k-2)!}{(m-s-i-1)!(m+n-i-r-1)!},
\end{eqnarray*}
where we have used the identity \citep[see][p. 509]{prudnikov:02}
\[
\sum_{k=0}^n (-1)^k {n \choose k}  {k+m \choose l}^{-1} = \frac{l}{n+l} {m+n \choose m-l}^{-1}.
\]
Substituting $S_{p0}$ and $S_{z0}$  in (\ref{joint3}) and by selecting $\eta=1$, we obtain
\begin{eqnarray*}
\lefteqn{P(A_{s}=k,B_{r}=i|H_0)=\frac{m!n!}{r!(n-k-r-1)!s!k!} }\\
& & \times \frac{(m+n-r-i-1)!(r+i)!}{i!(m+n)!}
\frac{(k+s)!(m+n-s-r-i-k-2)!}{(m-s-i-1)!(m+n-i-r-1)!} \\
& = &   \frac{\jbin{s+k}{s}\jbin{r+i}{r}}{\jbin{m+n}{n}} \jbin{m+n-s-r-i-k-2}{n-r-k-1}.
\end{eqnarray*}
Similarly, for $ m-s\le i \le m$ and $n-r \le k \le n$, the sums  $S_p'$ and $S_z'$ become
\begin{eqnarray*}
S_{p0}' & = & \frac{(m-i+n-r-1)!(k+s-m+i-n+r)!}{(k-n+r)!(k+s)!},
\\
  S_{z0}' &=& \frac{(k+s)!(m+n-s-k-1)!}{(n-k)!(m+n)!},
\end{eqnarray*}
and consequently, we get
\begin{eqnarray*}
\lefteqn{P(A_{s}=k,B_{r}=i|H_0)= \frac{m!n!}{(n-r-1)!(m-s-1)!(i-m+s)!(m-i)!} }\\
&  & \times\ \frac{(m-i+n-r-1)!(k+s-m+i-n+r)!(k+s)!(m+n-s-k-1)!}
{(k-n+r)!(k+s)!(n-k)!(m+n)!}
\\
& = &   \frac{\jbin{m+n-r-i-1}{n-r-1}\jbin{m+n-s-k-1}{m-s-1}}{\jbin{m+n}{n}}
 \jbin{k+i-m-n+s+r}{k-n+r}.
\end{eqnarray*}

The last case  of the theorem once again follows trivially due to (\ref{th-cas}).


\bigskip

\noindent{\bf Proof of Theorem \ref{large-appr}}

For fixed sample sizes $m$ and $n$, the lower tail of the exact distribution of $V_{\rho}$ is represented by $\sum_{i+k\le z} Q(k,i)$, where  $Q(k,i) = P(A_{s}=k,B_{r}=i|F=G)$.

As $m,n\to \infty$ and $m/n\to 1$, the behavior of
$\displaystyle
 \frac{ \jbin{m+n-s-r-i-k-2}{n-r-k-1}}{\jbin{m+n}{n}}
$
in $Q(k,i)$ is asymptotically equivalent to $2^{-(k+i+r+s+2)}$.  Therefore, the large-sample approximation of $\sum_{i+k\le z} Q(k,i)$ is given by
\begin{eqnarray*}
\sum_{i+k\le z} Q(k,i) & \sim &  \sum_{i+k\le z} \jbin{s+k}{s}\jbin{r+i}{r} 2^{-(k+i+r+s+2)}  \\
  && =  \sum_{i=0}^z \jbin{s+i}{s} 2^{-(i+s+1)} \sum_{k=0}^{z-i} \jbin{s+k}{s}2^{-(s+1)} 2^{-k}.
\end{eqnarray*}
The last sum represents the distribution function of a negative binomial random variable $\xi$ with parameters $s+1$ and 1/2. Using the well-known relationship between negative binomial distribution and binomial distribution \cite{bala-nev:03}, we have $P(\xi\le z) = P(\eta>s)$, where $\eta$ has binomial distribution with parameters $z+s+1$ and 1/2. Thus,

\begin{eqnarray*}
\lefteqn{\sum_{i+k\le z} Q(k,i)   \sim     \sum_{i=0}^z \jbin{s+i}{s} 2^{-(i+s+1)}\Big[1 -  \sum_{k=0}^{s} \jbin{z-i+s+1}{s}2^{-(z-i+s+1)} \Big]} \\
    && = \sum_{k=0}^z \jbin{z+s+1}{k} 2^{-(z+s+1)} - \sum_{i=0}^z \jbin{s+i}{s}
     \sum_{k=0}^s \jbin{z-i+s+1}{k} 2^{-(z+2s+2)} \\
    && =   -  \sum_{k=z+1}^{z+s+1} \jbin{z+2s+2}{k} 2^{-(z+2s+2)}+
     \sum_{i=0}^{z+s+1} \jbin{z+2s+2}{i} 2^{-(z+2s+2)}\\
      && =     \sum_{i=0}^{z} \jbin{z+2s+2}{i} 2^{-(z+2s+2)}
      =     \sum_{i=0}^{z} \jbin{2(s+1)+i-1}{i} 2^{-2(s+1)}2^{-i}.
\end{eqnarray*}

\end{document}